\documentclass[11pt]{article}
\title{On Karamata's proof of the Landau-Ingham Tauberian theorem}
\author{Michael M\"uger \\ Institute for Mathematics, Astrophysics and Particle Physics \\ 
Radboud University \\ Nijmegen, The Netherlands}

\usepackage{amssymb,amsmath,theorem,latexsym}

\newlength{\dinwidth}
\newlength{\dinmargin}
\def\1#1{{\bf #1}}
\def\2#1{{\cal #1}}
\def\3#1{{\sl #1}}
\def\4#1{{\tt #1}}
\def\5#1{{\sf #1}}
\def\6#1{{\mathfrak #1}}
\def\7#1{{\mathbb #1}}
\def\8#1{{\mathscr #1}}

\newcommand{\be}{\begin{equation}}
\newcommand{\ee}{\end{equation}}
\newcommand{\ba}{\begin{array}}
\newcommand{\ea}{\end{array}}
\newcommand{\bea}{\begin{eqnarray}}
\newcommand{\eea}{\end{eqnarray}}
 \newcommand{\bean}{\begin{eqnarray*}}
\newcommand{\eean}{\end{eqnarray*}}
\newcommand{\nn}{\nonumber}

\newcommand{\ve}{\varepsilon}

\newcommand{\impl}{\Rightarrow}
\newcommand{\rarr}{\rightarrow}

\newcommand{\DS}{\displaystyle}
\newcommand{\qed}{\hfill$\blacksquare$\\}

\def\endexem{\hfill{$\Box$}\medskip}

\theoremstyle{change}
\theoremheaderfont{\scshape}
\newtheorem{defin}{Definition}[section]
\newtheorem{defprop}{Definition/Proposition}[section]
\newtheorem{lemma}[defin]{Lemma}
\newtheorem{prop}[defin]{Proposition}
\newtheorem{theorem}[defin]{Theorem}
\newtheorem{coro}[defin]{Corollary}
\newtheorem{conj}[defin]{Conjecture}
\theorembodyfont{\rmfamily}
\newtheorem{remark}[defin]{Remark}

\newcommand{\bdefin}{\begin{defin}}
\newcommand{\blemma}{\begin{lemma}}
\newcommand{\bprop}{\begin{prop}}
\newcommand{\btheor}{\begin{theorem}}
\newcommand{\bcoro}{\begin{coro}}
\newcommand{\bdefprop}{\begin{defprop}}
\newcommand{\edefprop}{\end{defprop}}
\newcommand{\edefin}{\end{defin}}
\newcommand{\elemma}{\end{lemma}}
\newcommand{\eprop}{\end{prop}}
\newcommand{\etheor}{\end{theorem}}
\newcommand{\ecoro}{\end{coro}}
\newcommand{\bconj}{\begin{conj}}
\newcommand{\econj}{\end{conj}}
\newcommand{\brem}{\begin{remark}}
\newcommand{\erem}{\endexem\end{remark}}

\newcommand{\prf}{{\noindent\it Proof. }}
\def\mobj#1{\raise .4\unitlens\hbox{\put(0,0){$#1$}}}
\def\mychi{\raise 2pt\hbox{$\chi$}}

\numberwithin{equation}{section}

\begin{document}
\maketitle

\abstract{This is an exposition, in 12 pages including all prerequisites and a generalization, of
  Karamata's little known   elementary proof of the Landau-Ingham   Tauberian theorem, a result in
  real analysis from which   the Prime Number Theorem follows in a few   lines.}

\section{Introduction}
The aim of this paper is to give a self-contained, accessible and `elementary' proof of
of the following theorem, which we call the Landau-Ingham Tauberian theorem: 

\btheor \label{theor-LI} Let $f:[1,\infty)\rarr\7R$ be non-negative and non-decreasing and
assume that 
\be F(x):=\sum_{n\le x}f\left(\frac{x}{n}\right) \quad \mathrm{satisfies} \ \quad\
   F(x)=Ax\log x+Bx+C\frac{x}{\log x}+o\left(\frac{x}{\log x}\right). \label{eq-Kar}\ee
Then $f(x)=Ax+o(x)$, equivalently $f(x)\sim Ax$. 
\etheor

The interest of this theorem derives from the fact that, while ostensibly it is a result firmly
located in classical real analysis, the prime number theorem (PNT) $\pi(x)\sim\frac{x}{\log x}$
can be deduced from it by a few lines of Chebychev-style reasoning. (Cf.\ the Appendix.)

Versions of Theorem \ref{theor-LI} were proven by Landau \cite[\S 160]{landau} as early as 1909,
Ingham \cite[Theorem 1]{ingham}, Gordon \cite{gordon} and Ellison \cite[Theorem 3.1]{ellison}, but
none of these proofs was from scratch. Landau used as input the identity 
$\sum_n\frac{\mu(n)\log n}{n}=-1$. But the latter easily implies $M(x)=\sum_{n\le x}\mu(n)=o(x)$
which (as also shown by Landau) is equivalent to the PNT. Actually, 
$\sum_n\frac{\mu(n)\log n}{n}=-1$ is `stronger' than the PNT in the sense that it cannot be deduced
from the latter (other than by elementarily reproving the PNT with a sufficiently strong remainder
estimate). In this sense, Gordon's version of Theorem \ref{theor-LI} is an improvement, in that he
uses as input exactly the PNT (in the form $\psi(x)\sim x$) and thereby shows that Theorem
\ref{theor-LI} is not `stronger' than the PNT. Ellison's version assumes $M(x)=o(x)$
(and an $O(x^\beta)$ remainder with $\beta<1$ in (\ref{eq-Kar})). It is thus clear that none of
these approaches provides a proof of the PNT. Ingham's proof, on the other hand, departs from the
information that $\zeta(1+it)\ne 0$ (which can be deduced from the PNT, but also be proven {\it ab
  initio}). Thus his proof is not `elementary', but arguably it is one of the nicer and more
conceptual deductions of the PNT from $\zeta(1+it)\ne 0$ -- though certainly not the simplest (which
is \cite{zagier}) given that the proof requires Wiener's $L^1$-Tauberian theorem.

Our proof of Theorem \ref{theor-LI}  will essentially follow the elementary Selberg-style proof
given by Karamata 
\footnote{Note des \'editeurs : Jovan Karamata, n\'e pr\`es de Belgrade en 1902 et mort \`a Gen\`eve  en 1967, fut professeur \`a Gen\`eve d\`es 1951 et directeur de L'Enseignement Math\'ematique de   1954 \`a 1967. Voir M.\ Tomi\'c, \emph{Jovan Karamata (1902-1967)}, Enseignement  Math.\ \ (2)  {\bf 15}, 1--20 (1969).}
 \cite{karamata} under the assumption that $f$ is the summatory function of an
arithmetic function, i.e.\ constant between successive integers. We will remove this
assumption. For the proof of the PNT, this generality is not needed, but from an analysis
perspective it seems desirable, and it brings us fairly close to Ingham's version of the theorem,
which differed only in having $o(x)$ instead of $C\frac{x}{\log x}+o(\frac{x}{\log x})$ in the hypothesis.  

Unfortunately, Karamata's paper \cite{karamata} seems to be essentially forgotten: There are so few
references to it that we can discuss them all. It is mentioned in \cite{EI} by Erd\"os and Ingham
and in the book \cite{ellison} of Ellison and Mend\`es-France. (Considering that the latter authors
know Karamata's work, one may find it surprising that for their elementary proof of the PNT they
chose the somewhat roundabout route of giving a Selberg-style proof of $M(x)=o(x)$, using this to
prove a weak version of Theorem \ref{theor-LI}, from which then $\psi(x)\sim x$ is deduced.) Even
the two books \cite{BGT,korevaar} on Tauberian theory only briefly mention Karamata's
\cite{karamata} (or just the survey paper \cite{karamata2}) but then discuss in detail only Ingham's
proof. Finally, \cite{karamata,karamata2} are cited in the recent historical article \cite{nik}, but
its emphasis is on other matters. We close by noting that Karamata is not even mentioned in the only
other paper pursuing an elementary proof of a Landau-Ingham theorem, namely Balog's \cite{balog},
where a version of Theorem \ref{theor-LI} with a (fairly weak) error term in the conclusion is proven.  

Our reason for advertising Karamata's approach is that, in our view, it is the conceptually
cleanest and simplest of the Selberg-Erd\"os style proofs of the PNT, cf.\ \cite{selberg,erdos} and
followers, e.g.\ \cite{postn, nev, kalecki, levinson, schwarz, pollack}. For $f=\psi$ and
$f(x)=M(x)+\lfloor x\rfloor$, Theorem \ref{theor-LI} readily implies $\psi(x)=x+o(x)$ and
$M(x)=o(x)$. Making these substitutions in advance, the proof simplifies only marginally, but it
becomes less transparent (in particular for $f=\psi$) due to an abundance of non-linear
expressions. By contrast, Theorem \ref{theor-LI} is linear w.r.t.\ $f$ and $F$. To be sure, also the
proof given below has a non-linear core, cf.\ (\ref{eq-selb}) and Proposition \ref{prop-selb}, but
by putting the latter into evidence, the logic of the proof becomes clearer. One is actually led to
believe that the non-linear component of the proof is inevitable, as is also suggested by Theorem 2
in Erd\"os' \cite{erdos2}, to wit
\[ a_k\ge 0\ \forall k\ge 1\ \wedge\ \sum_{k=1}^N ka_k+\sum_{k+l\le   N}a_ka_l=N^2+O(1)\ \ \impl\ \
   \sum_{k=1}^N a_k=N+O(1),\]
from which the PNT can be deduced with little effort. (Cf.\ \cite{HT} for more in this direction.)

Another respect in which \cite{karamata} is superior to most of the later papers, including 
V.\ Nevanlinna's \cite{nev} (whose approach is adopted by several books \cite{schwarz,pollack}),
concerns the Tauberian deduction of the final result from a Selberg-style integral inequality. In
\cite{karamata}, this is achieved by a theorem attributed to Erd\"os (Theorem \ref{theor-2} below)
with clearly identified, obviously minimal hypotheses and an elegant proof. This advantage over
other approaches like \cite{nev}, which  tend to use further information about the discontinuities
of the function under consideration, is essential for our generalization to arbitrary non-decreasing
functions. However, we will have to adapt the proof (not least in order to work around an obscure issue).

In our exposition we make a point of avoiding the explicit summations over (pairs of) primes littering
many elementary proofs, almost obtaining a proof of the PNT free of primes! This is achieved by
defining the M\"obius and von Mangoldt functions $\mu$ and $\Lambda$  in terms of the functional
identities they satisfy and using their explicit computation only to show that they are bounded and
non-negative, respectively. Some of the proofs are  formulated in terms of parametric Stieltjes
integrals, typically of the form $\int f(x/t)dg(t)$ and integration by parts. We also do this in  
situations where $f$ and $g$ may both be discontinuous. Since our functions will always have 
bounded variation, thus at most countably many discontinuities, this can be justified by observing
that the resulting identities hold for all $x$ outside a countable set. Alternatively, we can
replace $f(x)$ at every point of discontinuity by $(f(x+0)+f(x-0))/2$ without changing the
asymptotics. For such functions, integration by parts always holds in the theory of
Lebesgue-Stieltjes integration, cf.\ \cite{hewitt,HS}.

The proof of Theorem \ref{theor-LI} exhibited below is, including all preliminaries, just 12 pages
long, and the  author hopes that this helps dispelling the prejudice that the elementary proofs of
the PNT are (conceptionally and/or technically) difficult.  Indeed he thinks that this is the most satisfactory of the 
elementary (and in fact of all) proofs of the PNT in that, besides not invoking complex analysis or
Riemann's $\zeta$-function, it minimizes number theoretic reasoning to a very well circumscribed
minimum. One may certainly dispute that this is desirable, but we will argue elsewhere that it is. 

The author is of course aware of the fact that the more direct elementary proofs of the PNT give
better control of the remainder term. (Cf.\ the review \cite{diamond} and the very recent paper
\cite{koukou}, which provides a ``a new and largely elementary proof of the best result known on the
counting function of primes in arithmetic progressions''.)  It is not clear whether this is necessarily so. \\

\noindent{\it Acknowledgments.} The author would like to thank the referees for constructive
comments that led to several improvements, in particular a better proof of Corollary \ref{coro-E2}.

\section{First steps and strategy}
\bprop \label{prop-estim}
Let $f:[1,\infty)\rarr\7R$ be non-negative and non-decreasing and assume that 
$F(x)=\sum_{n\le x}f(x/n)$ satisfies $F(x)=Ax\log x+Bx+o(x)$. Then 
\begin{itemize}
\item[(i)] $f(x)=O(x)$.
\item[(ii)] $\DS \int_{1-0}^x \frac{df(t)}{t}=A\log x+O(1)$.
\item[(iii)] $\DS\int_1^x\frac{f(t)-At}{t^2}dt=O(1)$.
\end{itemize}
\eprop

\prf (i) Following Ingham \cite{ingham}, we define $f$ to be $0$ on $[0,1)$ and compute
\bean f(x)-f\left(\frac{x}{2}\right)+f\left(\frac{x}{3}\right)-\cdots &=& F(x)-2F\left(\frac{x}{2}\right)\\
  &=& Ax\log x+Bx-2\left(A\frac{x}{2}\log\frac{x}{2}+B\frac{x}{2}\right)+o(x)=Ax\log 2+o(x).
\eean
With positivity and monotonicity of $f$, this gives $f(x)-f(x/2)\le  Kx$ for some $K>0$. Adding
these inequalities for $x, \frac{x}{2}, \frac{x}{4},\ldots$, we find $f(x)\le 2Kx$. Together with
$f\ge 0$, this gives (i). 

(ii) We compute
\bean F(x) &=& \sum_{n\le x} f\left(\frac{x}{n}\right)
   =\int_{1-0}^x f\left(\frac{x}{t}\right) d\lfloor t\rfloor \\
  &=& \left[\lfloor t\rfloor f\left(\frac{x}{t}\right)\right]_{t=1-0}^{t=x}
          -\int_{1-0}^x\lfloor t\rfloor\,df\left(\frac{x}{t}\right)\\
%  &=& \left(\lfloor x\rfloor f(1)-f(x)\right)
  &=& \lfloor x\rfloor f(1)
    -\int_{1-0}^x t\,df\left(\frac{x}{t}\right) +\int_{1-0}^x (t-\lfloor t\rfloor)\,df\left(\frac{x}{t}\right) \\
  &=& O(x) + \int_{1-0}^x \frac{x}{u}df(u) +\int_{1-0}^x (t-\lfloor t\rfloor)\,df\left(\frac{x}{t}\right). \\
\eean
In view of $0\le t-\lfloor t\rfloor<1$ and the weak monotonicity of $f$, the last integral is
bounded by $|\int_1^x df(x/t)|=f(x)-f(1)$, which is $O(x)$ by (i). Using the hypothesis about $F$, we
have 
\[ Ax\log x+Bx+o(x)=O(x)+x\int_{1-0}^x\frac{df(t)}{t}+O(x),\]
and division by $x$ proves the claim.

(iii) Integrating by parts, we have
\bean \int_1^x\frac{f(t)-At}{t^2}dt &=&
    -\frac{f(x)}{x}+\int_{1-0}^x \frac{df(t)}{t}-\int_1^x\frac{A}{t}dt \\
   &=& O(1)+(A\log x+O(1))-A\log x=O(1),\eean
where we used (i) and (ii).
\qed

\brem 1. The proposition can be proven under the weaker assumption $F(x)=Ax\log x+O(x)$, but we don't
bother since later we will need the stronger hypothesis anyway.

2. Theorem \ref{theor-LI}, which we ultimately want to prove, implies a strong form of
(iii): $\int_1^\infty\frac{f(t)-At}{t^2}dt=B-\gamma A$, cf.\ \cite{ingham}. Conversely, existence of
  the improper integral already implies $f(x)\sim Ax$, cf.\ \cite{zagier}.

3. Putting $f=\psi$ and using (\ref{eq-ch}), the above proofs of (i) and (ii) reduce to those of
Chebychev and Mertens, respectively. 
\erem

The following two theorems will be proven in Sections \ref{s-th1} and \ref{s-th2}, respectively:

\btheor \label{theor-1}
Let $f,F$ be as in Theorem \ref{theor-LI}. Then $g(x)=f(x)-Ax$ satisfies
\be \frac{|g(x)|}{x} \le \frac{1}{\log x}\int_1^x \frac{|g(t)|}{t^2}\,dt +o(1)\ \ \
   \mathrm{as}\ \ x\rarr\infty. \label{eq-ineq}\ee
(Here $f(x)\le g(x)+o(1)$ means that $f(x)\le g(x)+h(x)\ \forall x$, where $h(x)\rarr 0$ as $x\rarr\infty$.)
\etheor

\btheor \label{theor-2}
For $g:[1,\infty)\rarr\7R$, assume that there are $M, M'\ge 0$ such that 
\be x\mapsto g(x)+Mx\ \ \mbox{is non-decreasing}, \label{et1}\ee
\be \left| \int_1^x \frac{g(t)}{t^2} dt\right| \le M' \quad\forall x\ge 1. \label{et2}\ee
Then 
\be S:=\limsup_{x\rarr\infty}\frac{|g(x)|}{x}<\infty, \label{et3}\ee
and when $S>0$ we have
\be \limsup_{x\rarr\infty} \frac{1}{\log x}\int_1^x\frac{|g(t)|}{t^2}dt<S. \label{et4}\ee
\etheor

\brem \label{rem-th2} 1. Note that (\ref{et1}) implies that $g$ is Riemann integrable over finite
intervals. 

2. In our application, (\ref{et3}) already follows from Proposition
\ref{prop-estim} so that we do not need the corresponding part of the proof of Theorem
\ref{theor-2}. It will be proven nevertheless in order to give Theorem \ref{theor-2} an independent
existence. 
%
%3. With minimal additional effort, Theorem \ref{theor-2} can be proven under the  weaker
%assumption 
%\[ g(x')-g(x)\ge -M(x'-x)+o(x) \quad\forall 1\le x\le x'=O(x),\]
%cf.\ \cite{karamata}. We don't need this generality since our $g$ will satisfy (\ref{et1}) by
%construction. 
\erem

\noindent{\it Proof of Theorem \ref{theor-LI} assuming Theorems \ref{theor-1} and \ref{theor-2}.}
Since $f$ is nondecreasing, it is clear that $g(x)=f(x)-Ax$ satisfies (\ref{et1}) with $M=A$, and
(\ref{et2}) is implied by Proposition \ref{prop-estim}(iii). Now 
$S=\limsup|g(x)|/x$ is finite, by either Proposition \ref{prop-estim}(i) or the first conclusion of
Theorem \ref{theor-2}. Furthermore, $S>0$ would imply (\ref{et4}). But combining this with the
result (\ref{eq-ineq}) of Theorem \ref{theor-1}, we would have the absurdity 
\[ S=\limsup_{x\rarr\infty}\frac{|g(x)|}{x}\le \limsup_{x\rarr\infty}\frac{1}{\log x}\int_1^x \frac{|g(t)|}{t^2}\,dt<S.\]
Thus $S=0$ holds, which is equivalent to $\frac{g(x)}{x}=\frac{f(x)-Ax}{x}\rarr 0$,
as was to be proven.
\qed

The next two sections are dedicated to the proofs of Theorems \ref{theor-1} and \ref{theor-2}. The
statements of both results are free of number theory, and this is also the case for the proof of the
second. The proof of Theorem \ref{theor-1}, however, uses a very modest amount of number theory, but
nothing beyond M\"obius inversion and the divisibility theory of $\7N$ up to the fundamental theorem
of arithmetic.

\section{Proof of Theorem \ref{theor-1}}\label{s-th1}
\subsection{Arithmetic}\label{ss-mobius}
The aim of this subsection is to collect the basic arithmetic results that will be needed. We note
that this is very little.

We begin by noting that $(\7N,\cdot,1)$ is an abelian monoid. Given $n,m\in\7N$, we call $m$ a
divisor of $n$ if there is an $r\in\7N$ such that $mr=n$, in which case we  write $m|n$. In view of
the additive structure of the semiring $\7N$, it is clear that the monoid $\7N$ has cancellation 
($ab=ac\impl b=c$), so the quotient $r$ above is unique, and that the set of divisors 
of any $n$ is finite. 

Calling a function $f:\7N\rarr\7R$ an arithmetic function, the facts just stated allow us to define: 

\bdefin If $f,g:\7N\rarr\7R$ are arithmetic functions, their Dirichlet convolution $f\star g$
denotes the function 
\[ (f\star g)(n)=\sum_{d|n}f(d)g\left(\frac{n}{d}\right)=\sum_{a,b\atop ab=n} f(a)g(b).\]
\edefin

It is easy to see that Dirichlet convolution is commutative and associative. It has a unit given by the
function $\delta$ defined by $\delta(1)=1$ and $\delta(n)=0$ if $n\ne 1$. 

By $\11$ we denote the constant function $\11(n)=1$. Clearly, $(f\star\11)(n)=\sum_{d|n}f(d)$.

\blemma There is a unique arithmetic function $\mu$, called the M\"obius function, such that
$\mu\star\11=\delta$.  
\elemma

\prf $\mu$ must satisfy $\sum_{d|n}\mu(d)=\delta(n)$. Taking $n=1$ we see that $\mu(1)=1$. For
$n>1$ we have $\sum_{d|n}\mu(d)=0$, which is equivalent to 
\[ \mu(n)=-\sum_{d|n \atop d<n}\mu(d). \]
This uniquely determines $\mu(n)\in\7Z$ inductively in terms of $\mu(m)$ with $m<n$.
\qed

\bprop \label{prop-mu1} \begin{itemize}
\item[(i)] $\mu$ is multiplicative, i.e.\ $\mu(nm)=\mu(n)\mu(m)$ whenever $(n,m)=1$.
\item[(ii)] If $p$ is a prime then $\mu(p)=-1$, and $\mu(p^k)=0$ if $k\ge 2$.
\item[(iii)] $\mu(n)=O(1)$, i.e.\ $\mu$ is bounded. 
\end{itemize}
\eprop

\prf (i) Since $\mu(1)=1$, $\mu(nm)=\mu(n)\mu(m)$ clearly holds if $n=1$ or $m=1$. Assume, by way of
induction, that $\mu(uv)=\mu(u)\mu(v)$ holds whenever $(u,v)=1$ and $uv<nm$, and let $n\ne 1\ne m$
be relatively prime. Then every divisor of $nm$ is of the form $st$ with $s|n, t|m$, so that
\bean 0 &=& \sum_{d|nm} \mu(d)=\mu(nm)+\sum_{s|n, t|m\atop st<nm}\mu(st)
   =\mu(nm)+\sum_{s|n, t|m\atop st<nm}\mu(s)\mu(t) \\
   &=& \mu(nm)+\sum_{s|n}\mu(s)\sum_{t|m}\mu(t)-\mu(n)\mu(m)=\mu(nm)-\mu(n)\mu(m),
\eean
which is the inductive step. (ii) For $k\ge 1$, we have $\mu(p^k)=-\sum_{i=0}^{k-1}\mu(p^i)$,
inductively implying $\mu(p)=-1$ and $\mu(p^k)=0$ if $k\ge 2$. Thus $\mu(p^k)\in\{0,-1\}$, which
together with multiplicativity (i) gives $\mu(n)\in\{-1,0,1\}$ for all $n$, thus (iii).
\qed

\bprop \label{prop-Lambda1} 
\begin{itemize}
\item[(i)] The arithmetic function $\Lambda:=\log\star\mu$ is the unique solution of
  $\Lambda\star\11=\log$.  
\item[(ii)]  $\Lambda(n)=-\sum_{d|n}\mu(d)\log d$. In particular, $\Lambda(1)=0$.
\item[(iii)] $\Lambda(n)=\log p$ if $n=p^k$ where $p$ is prime and $k\ge 1$, and $\Lambda(n)=0$ otherwise. 
\item[(iv)] $\Lambda(n)\ge 0$.
\end{itemize}
\eprop

\prf (i) Existence: $\log\star\mu\star\11=\log\star\delta=\log$. Uniqueness: If
$\Lambda_1\star\11=\log=\Lambda_2\star\11$ then  
$\Lambda_1=\Lambda_1\star\delta=\Lambda_1\star\11\star\mu=\Lambda_2\star\11\star\mu=\Lambda_2\star\delta=\Lambda_2$. 

(ii) $\Lambda(n) =\sum_{d|n}\mu(d)\log\frac{n}{d}=\sum_{d|n}\mu(d)(\log n-\log d)
   =\log n\sum_{d|n}\mu(d)-\sum_{d|n}\mu(d)\log d$. Now use
   $\sum_{d|n}\mu(d)=\delta(n)$. $\Lambda(1)=0$ is obvious.

(iii) Using (ii), we have $\Lambda(p^k)=-\sum_{l=0}^k \mu(p^l) l\log p$, which together with 
Proposition \ref{prop-mu1}(ii) implies
$\Lambda(p^k)=\log p\ \forall k\ge 1$. If $n,m>1$ and  $(n,m)=1$ then by the    multiplicativity of $\mu$,
\bean \Lambda(nm) &=&-\sum_{s|n}\sum_{t|m}\mu(st)\log(st) =-\sum_{s|n}\sum_{t|m}\mu(s)\mu(t)(\log s+\log t)\\
 &=&\sum_{s|n}\mu(s)\log s\sum_{t|m}\mu(t)+\sum_{t|m}\mu(t)\log t\sum_{s|n}\mu(s)=0. \eean

(iv) Obvious consequence of (iii).
\qed

\brem The only properties of $\mu$ and $\Lambda$ needed for the proof of Theorem
\ref{theor-1} are the defining ones ($\mu\star\11=\delta,\ \Lambda\star\11=\log$), the trivial
consequence (ii) in Proposition \ref{prop-Lambda1}, and the boundedness of $\mu$ and the non-negativity
of $\Lambda$. 
In particular, the explicit computations of $\mu(n)$ and $\Lambda(n)$ in terms of the prime
factorization of $n$ were only needed to prove the latter two properties.
(Of course, the said properties of the functions $\mu$ and $\Lambda$ would be obvious if one
defined them by the explicit formulae proven above, but  this would be ad hoc and ugly, and one
would still need to use the fundamental theorem of arithmetic for proving that $\mu\star\11=\delta$
and $\Lambda\star\11=\log$.) 

Note that prime numbers will play no r\^ole whatsoever before we turn to the actual proof of the
prime number theorem in the Appendix,
%Section \ref{sec-pnt}, 
where the computation of $\Lambda(n)$ will be used
again.
\erem

\subsection{The (weighted) M\"obius transform}
\bdefin Given a function $f:[1,\infty)\rarr\7R$, its `M\"obius transform' is defined by
\[ F(x)=\sum_{n\le x}f\left(\frac{x}{n}\right).\]
\edefin

\blemma \label{lem-mobius}The M\"obius transform $f\mapsto F$ is invertible, the inverse M\"obius
transform being given by 
\[ f(x)=\sum_{n\le x} \mu(n)F\left(\frac{x}{n}\right).\]
\elemma

\prf We compute
\bean \sum_{n\le x}\mu(n) F\left(\frac{x}{n}\right) &=&
        \sum_{n\le x}\mu(n) \sum_{m\le x/n}f\left(\frac{x}{nm}\right)
  =\sum_{nm\le x}\mu(n)f\left(\frac{x}{nm}\right) \\
  &=& \sum_{r\le x} f\left(\frac{x}{r}\right) \sum_{s|r}\mu(s) 
        =\sum_{r\le x} f\left(\frac{x}{r}\right) \delta(r)=f(x),
\eean
where we used the defining property $\sum_{d|n}\mu(d)=\delta(n)$ of $\mu$.
\qed

\brem Since the point of Theorem \ref{theor-LI} is to deduce information about $f$ from information
concerning its M\"obius transform $F$, it is tempting to appeal to Lemma \ref{lem-mobius}
directly. However, in order for this to succeed, we would need control over 
$M(x)=\sum_{n\le  x}\mu(n)$, at least as good as $M(x)=o(x)$. But then one is back in Ellison's
approach mentioned in the introduction. The essential idea of the Selberg-Erd\"os approach to the
PNT, not entirely transparent in the early papers but clarified soon after \cite{TI}, is to consider
weighted M\"obius inversion formulae as follows.
\erem

\blemma \label{lem-TI1} Let $f:[1,\infty)\rarr\7R$ be arbitrary and $F(x)=\sum_{n\le x}f(x/n)$. Then 
\be f(x)\log x+\sum_{n\le x}\Lambda(n)f\left(\frac{x}{n}\right) 
   =\sum_{n\le  x}\mu(n)\log\frac{x}{n}\,F\left(\frac{x}{n}\right). \label{eq-TI}\ee
\elemma

\prf We compute
\[ \sum_{n\le x}\mu(n)\log\frac{x}{n}F\left(\frac{x}{n}\right) 
   =\log x \sum_{n\le x}\mu(n)F\left(\frac{x}{n}\right)-\sum_{n\le x}\mu(n)\log n\,F\left(\frac{x}{n}\right).\]
By Lemma \ref{lem-mobius}, the first term equals $f(x)\log x$, whereas for the second we have
\bean \sum_{n\le x}\mu(n)\log n\,F\left(\frac{x}{n}\right) &=& 
      \sum_{n\le x}\mu(n)\log n\,\sum_{m\le x/n} f\left(\frac{x}{nm}\right) 
   =\sum_{nm\le x}\mu(n)\log n\, f\left(\frac{x}{nm}\right) \\
   &=& \sum_{s\le x}\left( \sum_{n|s}\mu(n)\log n\right) \,f\left(\frac{x}{s}\right)
   =-\sum_{s\le x} \Lambda(s)  \,f\left(\frac{x}{s}\right),
\eean
the last equality being Proposition \ref{prop-Lambda1}(ii). Putting everything together,
we obtain (\ref{eq-TI}). 
\qed

\brem 1. Eq.\ (\ref{eq-TI}) is known as the `Tatuzawa-Iseki formula', cf.\ \cite[(8)]{TI} (and
\cite[p.\ 24]{karamata}). 

2. Without the factor $\log(x/n)$ on the right hand side, (\ref{eq-TI}) reduces to M\"obius
inversion. Thus (\ref{eq-TI}) is a sort of weighted M\"obius inversion formula. The presence 
of the sum involving $f(x/n)$ is very much wanted, since it will allow us to  obtain the 
integral inequality (\ref{eq-ineq}) involving all $f(t), t\in[1,x]$. In order to do so, we must get
rid of the explicit appearance of the function $\Lambda(n)$, which is very irregular and
about which we know little. This requires some preparations.
\erem

\blemma \label{lem-TI2} For any arithmetic function $f:\7N\rarr\7R$ we have
\[ f(n)\log n+\sum_{d|n}\Lambda(d)f\left(\frac{n}{d}\right) 
       =\sum_{d|n}\mu(d)\,\log\frac{n}{d}\,\sum_{m|(n/d)}f(m). \]
In particular, we have Selberg's identity:
\be \Lambda(n)\log n+\sum_{d|n}\Lambda(d)\Lambda\left(\frac{n}{d}\right)
      =\sum_{d|n}\mu(d)\log^2\frac{n}{d}. \label{eq-selb}\ee
\elemma

\prf If $f$ is an arithmetic function, i.e.\  defined only on $\7N$, we extend it to $\7R$ as being
  $0$ on $\7R\backslash\7N$. With this extension,
\[ F(n)=\sum_{m\le n}f\left(\frac{n}{m}\right) =\sum_{m|n}f\left(\frac{n}{m}\right) =\sum_{m|n}f(m),\]
so that (\ref{eq-TI}) becomes the claimed identity.  Taking $f(n)=\Lambda(n)$ and using
$\sum_{d|n}\Lambda(d)=\log n$, Selberg's formula follows. \qed

\subsection{Preliminary estimates}
\blemma The following elementary estimates hold as $x\rarr\infty$:
\bea
  \sum_{n\le x} \frac{1}{n} &=& \log x+\gamma+O\left(\frac{1}{x}\right), \label{s1}\\
  \sum_{n\le x} \frac{\log n}{n} &=& \frac{\log^2x}{2}+c+O\left(\frac{1+\log x}{x}\right), \label{s5}\\
  \sum_{n\le x} \log n &=& x\log x-x+O(\log x), \label{s2}\\
  \sum_{n\le x} \log\frac{x}{n} &=& x+O(\log x), \label{s2b}\\
  \sum_{n\le x} \log^2 n &=& x(\log^2x-2\log x+1)+ O(\log^2x),\label{s3}\\
 \sum_{n\le x} \log^2\frac{x}{n} &=& x+O(\log^2x). \label{s3b}
\eea
Here, $\gamma$ is Euler's constant and $c>0$.
\elemma

\prf (\ref{s1}): We have
\[ \sum_{n=1}^N\frac{1}{n}-\int_1^N\frac{dt}{t}=\int_{1-0}^N \frac{d(\lfloor t\rfloor -t)}{t}
   =\left[\frac{\lfloor t\rfloor -t}{t}\right]_1^N+\int_1^N \frac{t-\lfloor t\rfloor}{t^2}dt. \]
Since $0\le t-\lfloor t\rfloor<1$, the integral on the r.h.s.\ converges as $N\rarr\infty$ to some
number $\gamma$ (Euler's constant) strictly between 0 and $1=\int_1^\infty dt/t^2$. Thus
\[ \sum_{n=1}^N\frac{1}{n}=\int_1^N\frac{dt}{t}+\gamma-\int_N^\infty \frac{t-\lfloor t\rfloor}{t^2}dt 
  =\log N+\gamma+O\left(\frac{1}{N}\right). \]

(\ref{s5}): Similarly to the proof of (\ref{s1}), we have
\[ \sum_{n=1}^N\frac{\log n}{n}-\int_1^N\frac{\log t}{t}dt=\int_{1-0}^N \frac{\log t}{t}\,d(\lfloor t\rfloor -t)
   =\left[\frac{(\lfloor t\rfloor -t)\log t}{t}\right]_1^N+\int_1^N \frac{(t-\lfloor t\rfloor)\log t}{t^2}dt. \]
The final integral converges to some $c>0$ as $N\rarr\infty$ since $(\log t)/t^2=O(t^{-2+\ve})$. Using
\[ \int_1^x \frac{\log t}{t}dt= \frac{\log^2 x}{2},\quad\quad
  \int_N^\infty\frac{\log t}{t^2}dt=-\left[\frac{\log   t}{t}\right]_N^\infty+\int_N^\infty\frac{dt}{t^2} =\frac{1+\log N}{N} \]
we have 
\[ \sum_{n=1}^N\frac{\log n}{n}=\int_1^x \frac{\log t}{t}dt+c-\int_N^\infty \frac{(t-\lfloor t\rfloor)\log t}{t^2}dt
   =\frac{\log^2 x}{2} + c + O\left(\frac{1+\log x}{x}\right). \]

(\ref{s2}): By monotonicity, we have
\[ \int_1^x \log t\,dt \le \sum_{n\le x} \log n \le \int_1^{x+1} \log t\,dt. \]
Combining this with $\int_1^x\log t\,dt=x\log x-x+1$, (\ref{s2}) follows.

(\ref{s2b}): Using (\ref{s2}), we have
\[ \sum_{n\le x} \log\frac{x}{n}= \lfloor x\rfloor\log x-\sum_{n\le x} \log n
   = (x+O(1)) \log x-(x\log x-x+O(\log x))=x+O(\log x). \]

(\ref{s3}): By monotonicity, 
\[ \int_1^x \log^2t\,dt \le \sum_{n\le x} \log^2 n \le \int_1^{x+1} \log^2t\,dt. \]
Now, 
\[ \int_1^x\log^2t\,dt=\int_0^{\log x} e^u u^2du=[e^u(u^2-2u+1)]_0^{\log x}=x(\log^2x-2\log x+1)-1.\]
Combining these two facts, (\ref{s3}) follows.

 (\ref{s3b}): Using (\ref{s2}) and (\ref{s3}), we compute
\begin{flalign*}
\quad\sum_{n\le x} \log^2\frac{x}{n} &= \sum_{n\le x} (\log x-\log n)^2 & \\
\quad  &= \lfloor x\rfloor\log^2x-2\log x(x\log x-x+O(\log x))+x(\log^2x-2\log x+1)+ O(\log^2x) & \\
\quad  &= x+O(\log^2x) & \blacksquare
\quad\end{flalign*}

\vspace{.2cm}

\bprop \label{prop-mu2}
The following estimates involving the M\"obius function hold as $x\rarr\infty$:
\bea
 \sum_{n\le x} \frac{\mu(n)}{n} &=&  O(1), \label{mu1}\\
 \sum_{n\le x} \frac{\mu(n)}{n}\log\frac{x}{n} &=& O(1), \label{mu2}\\
 \sum_{n\le x} \frac{\mu(n)}{n}\log^2\frac{x}{n} &=& 2\log x+ O(1). \label{mu3}
\eea
\eprop

\prf (\ref{mu1}): If $f(x)=1$ then $F(x)=\lfloor x\rfloor$. M\"obius inversion (Lemma
\ref{lem-mobius}) gives
\be 1=\sum_{n\le x}\mu(n)\left\lfloor\frac{x}{n}\right\rfloor=\sum_{n\le  x}\mu(n)\left(\frac{x}{n}+O(1)\right)
   =x\sum_{n\le  x}\frac{\mu(n)}{n}+\sum_{n\le  x}O(1),\label{eq-mu}\ee
where we used $\mu(n)=O(1)$ (Proposition \ref{prop-mu1}(iii)). In view of $\sum_{n\le x}O(1)=O(x)$, we have 
$\sum_{n\le x}\mu(n)/n=O(x)/x=O(1)$.  

(\ref{mu2}): If $f(x)=x$ then $F(x)=\sum_{n\le x}x/n=x\log x+\gamma x+O(1)$ by (\ref{s1}). By
M\"obius inversion,
\[ x=\sum_{n\le x}\mu(n)\left( \frac{x}{n}\log \frac{x}{n}+\gamma \frac{x}{n}+O(1)\right)
   =x \sum_{n\le x}\frac{\mu(n)}{n}\log \frac{x}{n}+xO(1)+O(x), \]
where we used (\ref{mu1}) and Proposition \ref{prop-mu1}(iii). From this we easily read off (\ref{mu2}).

(\ref{mu3}): If $f(x)=x\log x$ then 
\bean  F(x) &=& \sum_{n\le x} \frac{x}{n}\log \frac{x}{n}=\sum_{n\le x} \frac{x}{n}(\log x-\log n)
   =x\log x\sum_{n\le x}\frac{1}{n}-x\sum_{n\le x}\frac{\log n}{n} \\
  &=& x\log x\left(\log x+\gamma+O\left(\frac{1}{x}\right)\right)
     -x\left(\frac{\log^2x}{2}+c+O\left(\frac{1+\log x}{x}\right)\right) \\
  &=& \frac{1}{2}x\log^2x+\gamma x\log x-cx+O(1+\log x),
\eean
by (\ref{s1}) and (\ref{s5}). Now M\"obius inversion gives
\bean x\log x &=& \sum_{n\le x}\mu(n)\left(\frac{x}{2n}\log^2\frac{x}{n}    
      +\gamma\frac{x}{n}\log \frac{x}{n}-c\frac{x}{n}+O(1+\log \frac{x}{n})\right)\\
   &=& \frac{x}{2}\sum_{n\le x}\frac{\mu(n)}{n}\log^2\frac{x}{n} +xO(1)+xO(1)+O(x),   
\eean
where we used (\ref{mu1}), (\ref{mu2}) and (\ref{s2b}), and division by $x/2$ gives (\ref{mu3}).
\qed

\subsection{Conclusion}
\bprop [Selberg, Erd\"os-Karamata \cite{EK}] \label{prop-selb} Defining
\[ K(1)=0,\quad\quad K(n)=\frac{1}{\log n}\sum_{d|n}\Lambda(d)\Lambda\left(\frac{n}{d}\right)
\ \ \mathrm{if}\ \ n\ge 2,\label{eq-K}\]
we have $K(n)\ge 0$ and
\be \sum_{n\le x} (\Lambda(n)+K(n))=2x+O\left(\frac{x}{\log x}\right). \label{eq-K2}\ee
\eprop

\prf The first claim is obvious in view of Proposition \ref{prop-Lambda1}(iv). We estimate
\bean U(x) &:= &\sum_{n\le x}\sum_{d|n}\mu(d)\log^2\frac{n}{d}
   \ =\ \sum_{n\le x}\mu(n) \sum_{m\le x/n}\log^2 m \\
  &=& \sum_{n\le x}\mu(n) \left( \frac{x}{n}(\log^2\frac{x}{n}-2\log\frac{x}{n}+1)+ O(\log^2\frac{x}{n})\right)\\
  &=& x(2\log x+O(1))-2x O(1)+ xO(1)+ O(x)=2x\log x+O(x). 
\eean
Here we used (\ref{s3}), (\ref{mu3}), (\ref{mu2}), (\ref{mu1}), the fact $\mu(n)=O(1)$, and (\ref{s3b}). 
Comparing (\ref{eq-K2}) and (\ref{eq-selb}), we have
\bean \sum_{n\le x} \Lambda(n)+K(n) &=& \sum_{2\le n\le x}\frac{1}{\log n} \sum_{d|n} \mu(d)\log^2\frac{n}{d}\\
   &=& \int_{2-0}^x\frac{dU(t)}{\log t}\ =\ \left[\frac{U(t)}{\log t}\right]_2^x +\int_2^x \frac{U(t)}{t\log^2t}dt \\
   &=& 2x+O\left(\frac{x}{\log x}\right)+ \int_2^x\frac{dt}{\log t}+O\left(\int_2^x\frac{dt}{\log^2t}\right).
\eean 
In view of the estimate
\[ \int_2^x\frac{dt}{\log t}=\int_2^{\sqrt{x}}\frac{dt}{\log t}+\int_{\sqrt{x}}^x\frac{dt}{\log t}
    \le \frac{\sqrt{x}}{\log 2}+\frac{x}{\log\sqrt{x}}=O\left(\frac{x}{\log x}\right), \]
we are done.
\qed

\brem In view of (\ref{eq-selb}), the above estimate $U(x)=2x\log x+O(x)$ is equivalent to
\[ \sum_{n\le x} \Lambda(n)\log n+\sum_{ab\le x}\Lambda(a)\Lambda(b)=2x\log x+O(x),\]
which is used in most Selberg-style proofs. (It would lead to (\ref{eq-balog}) with $k=2$.)
\erem

\bprop \label{prop-G} If $g:[1,\infty)\rarr\7R$ is such that
\be G(x)=\sum_{n\le x}g\left(\frac{x}{n}\right)
        =Bx+C\frac{x}{\log x}+o\left(\frac{x}{\log x}\right) \label{eq-Kar2}\ee
then 
\be g(x)\log x+\sum_{n\le x}\Lambda(n)\,g\left(\frac{x}{n}\right)= o(x\log x). \label{e2}\ee
\eprop

\prf In view of Lemma \ref{lem-TI1}, all we have to do is estimate 
\[ \sum_{n\le  x}\mu(n)\log\frac{x}{n}\left( B\frac{x}{n}+C\frac{x}{n\log \frac{x}{n}}
        +o\left(\frac{x}{n\log \frac{x}{n}}\right) \right)=S_1+S_2+S_3. \]
The three terms are
\bean S_1 &=& Bx \sum_{n\le  x}\frac{\mu(n)}{n}\log\frac{x}{n}=x O(1)=O(x), \\
   S_2 &=& Cx \sum_{n\le  x}\frac{\mu(n)}{n}=xO(1)=O(x), \\
   S_3 &=& \sum_{n\le  x}\mu(n) o\left(\frac{x}{n}\right) =\sum_{n\le  x}o\left(\frac{x}{n}\right)
   =o\left(x\sum_{n\le  x}\frac{1}{n}\right)=o(x\log x), 
\eean
where we used (\ref{mu2}), (\ref{mu1}), and $\mu(n)=O(1)$, respectively.
\qed

\noindent{\it Proof of Theorem \ref{theor-1}.} 
In view of $g(x)=f(x)-Ax$ and Proposition \ref{prop-estim} (i), (ii), we immediately have
\be g(x)=O(x), \quad\quad\quad \int_1^x\frac{dg(u)}{u}=O(1). \label{eq-g}\ee
Furthermore, since $f$ satisfies (\ref{eq-Kar}), and (\ref{s1}) gives 
$\sum_{n\le x} Ax/n=Ax\log x+A\gamma x+O(1)$, the M\"obius transform $G$ of $g(x)=f(x)-Ax$ satisfies 
(\ref{eq-Kar2}) (with a different $B$), so that Proposition \ref{prop-G} applies and  (\ref{e2}) holds. 

Writing 
$N(x)=\sum_{n\le   x}\Lambda(n)+K(n)$, by Proposition \ref{prop-selb} we have $N(x)=2x+\omega(x)$
with $\omega(x)=o(x)$. Now,
\bea \sum_{n\le  x}(\Lambda(n)+K(n))g\left(\frac{x}{n}\right)
    &=& \int_{1-0}^xg\left(\frac{x}{t}\right)dN(t) \nn\\
  &=& \left[ N(t)g\left(\frac{x}{t}\right)\right]_{1-0}^x-\int_1^x N(t)dg\left(\frac{x}{t}\right)\nn\\
  &=& (N(x)g(1)-N(1-0)g(x)) +\int_1^x N\left(\frac{x}{u}\right) dg(u)\nn\\
  &=& O(x) +\int_1^x\left(\frac{2x}{u}+o\left(\frac{2x}{u}\right)\right)dg(u) \nn\\
  &=& O(x) +2x\int_1^x \frac{dg(u)}{u}+o\left(x\int_1^x \frac{dg(u)}{u}\right)\nn\\
  &=& O(x)+O(x)+o(x)=O(x), \label{eq-x1}
\eea
where we used  (\ref{eq-g}).
%used $g(x)=O(x)$ and $\int_1^x\frac{dg(u)}{u}=O(1)$, both of which are immediate
%consequences of $g(x)=f(x)-Ax$ and Proposition \ref{prop-estim}(i) and (ii). 
On the other hand, 
\bea  \lefteqn{ \sum_{n\le x} (\Lambda(n)+K(n))\,\left|g\left(\frac{x}{n}\right)\right| = 
    \int_{1-0}^x \left|g\left(\frac{x}{t}\right)\right| \, dN(t) } \nn\\
   &=& 2\int_1^x \left|g\left(\frac{x}{t}\right)\right|\,dt+\int_{1-0}^x \left|g\left(\frac{x}{t}\right)\right|\,d\omega(t) \nn\\
   &=& 2x\int_1^x \frac{|g(t)|}{t^2}dt -\int_{1-0}^x \omega(t)\, d\!\left|g\left(\frac{x}{t}\right)\right| +\left[g\left(\frac{x}{t}\right)\omega(t)\right]_{t=1-0}^{t=x}\nn\\
   &=& 2x\int_1^x \frac{|g(t)|}{t^2}dt +\int_1^{x+0} \omega\left(\frac{x}{t}\right)\, d|g(t)| +g(1)\omega(x)-g(x+0)\omega(1-0).\label{eq-x2}
\eea
In view of $g(x)=O(x)$ and $\omega(x)=o(x)$, the sum of the last two terms is $O(x)$. Furthermore,
\bean \int_1^{x+0} \omega\left(\frac{x}{t}\right)\,d|g(t)| &=& o\left( x\int_1^{x+0}\frac{|d|g(t)||}{t}\right)
   \le o\left( x\int_1^{x+0}\frac{|dg(t)|}{t}\right) \\
   &\le& o\left( x\int_1^{x+0}\frac{df+Adt}{t}\right)=o(x\log x), \eean
where we used $g(x)=f(x)-Ax$ and $df=|df|$ (since $f$ is non-decreasing) to obtain 
$|dg|=|df-Adt|\le |df|+Adt=df+Adt$ and Proposition \ref{prop-estim}(ii). Plugging this into
(\ref{eq-x2}), we have 
\be \sum_{n\le x} (\Lambda(n)+K(n))\,\left|g\left(\frac{x}{n}\right)\right| =
       2x\int_1^x \frac{|g(t)|}{t^2}dt+ o(x\log x). \label{eq-x3}\ee
After these preparations, we can conclude quickly: Subtracting (\ref{eq-x1}) from (\ref{e2}) we obtain
\[ g(x)\log x= \sum_{n\le x} K(n) g\left(\frac{x}{n}\right) + o(x\log x).\]
Taking absolute values of this and of (\ref{e2}) while observing that $\Lambda$ and $K$ are
non-negative, we have the inequalities
\[ |g(x)|\log x\le \sum_{n\le x} \Lambda(n) \left|g\left(\frac{x}{n}\right)\right| + o(x\log x), \quad\quad
    |g(x)|\log x\le \sum_{n\le x} K(n) \left|g\left(\frac{x}{n}\right)\right| + o(x\log x).\]
Adding these inequalities and comparing with (\ref{eq-x3}) we have
\[ 2|g(x)|\log x \le \sum_{n\le x}(\Lambda(n)+K(n)) \left|g\left(\frac{x}{n}\right)\right| + o(x\log x)
    = 2x\int_1^x \frac{|g(t)|}{t^2}dt+ o(x\log x), \]
so that (\ref{eq-ineq}), and with it Theorem \ref{theor-1}, is obtained dividing by $2x\log x$.
\qed

\brem 1. We did not use the full strength of Proposition \ref{prop-selb}, but only an
$o(x)$ remainder.

2. Inequality (\ref{eq-ineq}) is the special case $k=1$ of the more general integral inequality
\be \frac{|g(x)|}{x}\log^kx\le k\int_1^x\frac{|g(t)|\log^{k-1}t}{t^2}dt+O(\log^{k-c}x)
    \quad\forall k\in\7N \label{eq-balog}\ee
proven in \cite{balog}, assuming a $O\left(\frac{x}{\log^2x}\right)$ in (\ref{eq-Kar}) instead of
$C\frac{x}{\log x}+o\left(\frac{x}{\log x}\right)$. 
\erem

\section{Proof of Theorem \ref{theor-2}}\label{s-th2}
The proof will be based on the following proposition, to be proven later:

\bprop \label{prop-E} If $s:[0,\infty)\rarr\7R$ satisfies
\be e^{t'}s(t')-e^t s(t)\ge -M(e^{t'}-e^t) \quad \forall t'\ge t\ge 0,\label{eq-s1}\ee
\be \left|\int_0^x s(t)dt\right|\le M'\quad\forall x\ge 0, \label{eq-s2}\ee
and   $S=\lim\sup|s(x)|>0$ then there exist numbers $0<S_1<S$ and $e,h>0$ such that  
\be \mu(E_{x,h,S_1}) \ge e\quad \forall x\ge 0,\quad \mathrm{where}\quad
  E_{x,h,S_1}=\{ t\in [ x,x+h]\ | \ |s(t)|\le S_1 \},\label{eq-E}\ee
and $\mu$ denotes the Lebesgue measure.
\eprop

\noindent{\it Proof of Theorem \ref{theor-2} assuming Proposition \ref{prop-E}.}
It is convenient to replace $g:[1,\infty)\rarr\7R$ by $s:[0,\infty)\rarr\7R,\ s(t)= e^{-t} g(e^t)$.
Now $s$ is locally integrable, and the assumptions (\ref{et1}) and (\ref{et2}) become
%\be e^{t'}s(t')-e^t s(t)\ge -M(e^{t'}-e^t)+o(e^t)\quad \forall 0\le t\le t'=t+O(t),\label{eq-s1}\ee
(\ref{eq-s1}) and (\ref{eq-s2}), respectively, 
%\be e^{t'}s(t')-e^t s(t)\ge -M(e^{t'}-e^t) \quad \forall t'\ge t\ge 0,\label{eq-s1}\ee
%\be \left|\int_0^x s(t)dt\right|\le M'\quad\forall x\ge 0, \label{eq-s2}\ee
whereas the conclusions (\ref{et3}) and (\ref{et4}) assume the form
\be S=\limsup_{t\rarr\infty}|s(t)|<\infty, \label{eq-s3}\ee
\be S>0 \quad\impl\quad \limsup_{x\rarr\infty} \frac{1}{x}\int_0^x |s(t)|dt < S. \label{eq-s4}\ee
The proof of (\ref{eq-s3}) is easy: Dividing (\ref{eq-s1}) by $e^{t'}$ and integrating over
$t'\in[t,t+h]$, where $h>0$, one obtains
%\[ \int_t^{t+h}  s(t')dt' - s(t)(1-e^{-h})\ge -Mh+M(1-e^{-h})+o(1), \] 
\[ \int_t^{t+h}  s(t')dt' - s(t)(1-e^{-h})\ge -Mh+M(1-e^{-h}), \] 
and using $|\int_a^b s(t)dt|\le |\int_0^as(t)dt|+|\int_0^bs(t)dt|\le 2M'$ by (\ref{eq-s2}), we have
the upper bound
%\[ s(t)\le \frac{2M'+M(e^{-h}-1+h)}{1-e^{-h}} +o(1). \] 
\[ s(t)\le \frac{2M'+M(e^{-h}-1+h)}{1-e^{-h}} . \] 
Similarly, dividing (\ref{eq-s1}) by $e^t$ and integrating over $t\in[t'-h,t']$, one obtains the
lower bound
%\[ -\frac{2M'+M(e^h-1-h)}{e^h-1} +o(1)\le s(t),\] 
\[ -\frac{2M'+M(e^h-1-h)}{e^h-1} \le s(t),\] 
thus (\ref{eq-s3}) holds. 

Assuming $S>0$, let $S_1,h,e$ as provided by Proposition \ref{prop-E}. For each $\widehat{S}>S$
there is $x_0$ such that $x\ge x_0\impl |s(x)|\le\widehat{S}$. Given $x\ge x_0$ and putting
$N=\left\lfloor\frac{x-x_0}{h}\right\rfloor$, we have  
\bean \int_0^x |s(t)|dt &=& \int_0^{x_0}|s(t)|dt+ \sum_{n=1}^N \int_{x_0+(n-1)h}^{x_0+nh} |s(t)|dt
     +\int_{x_0+Nh}^x |s(t)|dt  \\
   &\le & 2M'+ N [ \widehat{S}(h-e)+S_1e] + 2M' \\
   &=& \left(\frac{x-x_0}{h}+O(1)\right) h \left[\left(1-\frac{e}{h}\right)\widehat{S}+\frac{e}{h}S_1\right] +4M'\\
  &=& x \left[\left(1-\frac{e}{h}\right)\widehat{S}+\frac{e}{h}S_1\right]   +O(1).
\eean
Thus
\[ \limsup_{x\rarr\infty}\frac{1}{x}\int_0^x |s(t)|dt\le \left(1-\frac{e}{h}\right)\widehat{S}+\frac{e}{h}S_1.  \]
Since $S_1<S$ and since $\widehat{S}>S$ can be chosen arbitrarily close to $S$, (\ref{eq-s4}) holds
and thus Theorem \ref{theor-2}.
\qed

In order to make plain how the assumptions (\ref{eq-s1}) and (\ref{eq-s2}) enter the proof of 
Proposition \ref{prop-E}, we prove two intermediate results that each use only one of the
assumptions. For the first we need a ``geometrically obvious'' lemma of isoperimetric character:

\blemma Let $t_1<t_2,\ C_1>C_2>0$ and $k:[t_1,t_2]\rarr\7R$ non-decreasing with
$k(t_1)\ge C_1e^{t_1}$ and $k(t_2)\le C_2e^{t_2}$. Then 
\[ \mu\left(\{ t\in[t_1,t_2] \ | \ C_2e^t\le k(t)\le C_1 e^t\}\right) \ge \log\frac{C_1}{C_2}.\]
\elemma

\prf As a non-decreasing function, $k$ has left and right limits $k(t\pm 0)$ everywhere and
$k(t-0)\le k(t)\le k(t+0)$. The assumptions imply 
$t_1\in A:=\{ t\in[t_1,t_2] \ | \ k(t)\ge C_1e^t\}$, thus we can define $T_1=\sup(A)$. 
Quite obviously we have $t>T_1\impl k(t)<C_1e^t$, which together with the non-decreasing property of
$k$ and the continuity of the exponential function implies $k(T_1+0)\le C_1e^{T_1}$ (provided
$T_1<t_2$). We have $T_1\in A$ if and only if $k(T_1)\ge C_1e^{T_1}$. If $T_1\not\in A$ then
$T_1>t_1$, and every  interval $(T_1-\ve,T_1)$ (with $0<\ve<T_1-t_1$) contains points $t$ such that
$k(t)\ge C_1e^t$. This implies  $k(T_1-0)\ge C_1e^{T_1}$. Now assume $T_1=t_2$. If $T_1\in A$ then
$C_1e^{T_1}\le k(T_1)\le C_2e^{T_1}$. If $T_1\not\in A$ then 
$C_1e^{T_1}\le k(T_1-0)\le k(T_1)\le C_2e^{t_2}$. In both cases we arrive at a contradiction since $C_2<C_1$.
Thus $T_1<t_2$. 
If $T_1\in A$ (in particular if $T_1=t_1$) then 
$C_1e^{T_1}\le k(T_1)\le k(T_1+0)\le C_1e^{T_1}$. Thus $k$ is continuous from the right at $T_1$ and  
$k(T_1)=C_1e^{T_1}$. If $T_1\not\in A$ then $T_1>t_1$ and
$C_1e^{T_1}\le k(T_1-0)\le k(T_1+0)\le C_1e^{T_1}$. This implies $k(T_1)=C_1e^{T_1}$, thus the
contradiction $T_1\in A$. Thus we always have $T_1\in A$, thus $k(T_1)=C_1e^{T_1}$.

Now let $B=\{ t\in[T_1,t_2]\ | \  k(t)\le C_2e^t\}$. We have $t_2\in B$, thus $T_2=\inf(B)$ is 
defined and $T_2\ge T_1$. Arguing similarly as before we have $t<T_2\impl k(t)>C_2e^t$,   implying
$k(T_2-0)\ge C_2e^{T_2}$. And if $T_2<t_2$ and $T_2\not\in B$ then $k(T_2+0)\le C_2e^{T_2}$.
If $T_2\in B$ (in particular if $T_2=t_2$) then  $C_2e^{T_2}\le k(T_2-0)\le k(T_2)\le  C_2e^{T_2}$,
implying $k(T_2-0)=k(T_2)=C_2e^{T_2}$ so that $k$ is continuous from the left at $T_2$. If
$T_2\not\in B$ then $T_2<t_2$ and $C_2e^{T_2}\le k(T_2-0)\le k(T_2+0)\le C_2e^{T_2}$, implying
$k(T_2)=C_2e^{T_2}$ and thus a contradiction. Thus we always have $T_2\in B$, thus $k(T_2)=C_2e^{T_2}$. 

%If $T_1=T_2$ then $C_1e^{T_1}=C_2e^{T_1}$, which is impossible since $C_2<C_1$. Thus $T_1<T_2$.  

By the above results, we have $C_2e^t\le k(t)\le C_1 e^t\ \forall t\in[T_1,T_2]$ and thus
\be \mu\left(\{   t\in[t_1,t_2] \ | \ C_2e^t\le k(t)\le C_1 e^t\}\right) \ge T_2-T_1. \label{mu}\ee
Using once  more that $k$ is non-decreasing, we have 
\[ C_1e^{T_1}= k(T_1)\le k(T_2)=C_2e^{T_2}, \]
implying $T_2-T_1\ge \log\frac{C_1}{C_2}$, and combining this with (\ref{mu}) proves the claim.
\qed

\bcoro \label{coro-E2} Assume that $s:[0,\infty)\rarr\7R$ satisfies (\ref{eq-s1}) and  
$s(t_1)\ge S_1\ge S_2\ge s(t_2)$, where $S_2+M>0$. Then
%\[ \mu([t_1,t_2]\cap s^{-1}([S_2,S_1]))\ge\log\frac{S_1+M}{S_2+M}. \]
\[ \mu(\{t\in [t_1,t_2]\ | \ s(t)\in [S_2,S_1] \})   \ge\log\frac{S_1+M}{S_2+M}. \]
\ecoro

\prf We note that (\ref{eq-s1}) is equivalent to the statement that the function 
$k:t\mapsto e^t(s(t)+M)$ is non-decreasing. The assumption  $s(t_1)\ge S_1\ge S_2\ge s(t_2)$ implies 
$k(t_1)\ge (S_1+M)e^{t_1}$ and $k(t_2)\le(S_2+M)e^{t_2}$.  Now the claim follows directly by an
application of the preceding lemma.
\qed

\blemma \label{lem-E1} Let $s:[0,\infty)\rarr\7R$ be integrable over bounded intervals, satisfying
(\ref{eq-s2}). Let $e>0$ and $0<S_2<S_1$ be arbitrary, and assume
\be h\ge 2\left( e+\frac{M'}{S_1}+\frac{M'}{S_2}\right).\label{eq-h}\ee
Then every interval $[x,x+h]$ satisfies at least one of the following conditions:
\begin{itemize}
\item[(i)] $\DS\mu(E_{x,h,S_1}) \ge e$, where $E_{x,h,S_1}$ is as in (\ref{eq-E}),
\item[(ii)] there exist $t_1,t_2$ such that $x\le t_1<t_2\le x+h$ and $s(t_1)\ge S_1$ and 
$s(t_2)\le S_2$.
\end{itemize}
\elemma

\prf It is enough to show that falsity of (i) implies (ii). Define
\[ T=\sup\{ t\in[x,x+h]\ | \ s(t)\le S_2\}, \]
with the understanding that $T=x$ if $s(t)>S_2$ for all $t\in[x,x+h]$.
Then $s(t)>S_2\ \forall t\in(T,x+h]$, which implies
\[  (x+h-T)S_2\le  \int_{T}^{x+h} s(t)dt\le 2M' \]
and therefore
\be x+h-T\ \le \frac{2M'}{S_2}. \label{eq-t2}\ee
%\be x+h-T\ \le\ 2M'/S_2. \label{eq-t2}\ee
We observe that (\ref{eq-t2}) with $T=x$  would contradict (\ref{eq-h}). Thus
$x<T\le x+h$, so we can indeed find a $t_2\in[x,x+h]$ with $s(t_2)\le S_2$. Since we do not assume
continuity of $s$, we cannot claim that we may take $t_2=T$, but by definition a $t_2$ can be found in 
$(T-\ve,T]$ for every $\ve>0$.

Now we claim that there is a point $t_1\in[x,t_2]$ such that $s(t_1)\ge S_1$. Otherwise, we would have
$s(t)<S_1$ for all $t\in[x,t_2]$. By definition, $|s(t)|\le S_1$ for $t\in E_{x,h,S_1}$, thus
$|s|>S_1$ on the complement of $E_{x,h,S_1}$. Combined with $s(t)<S_1$ for $t\in[x,t_2]$, this
means $s(t)<-S_1$ whenever $t\in[x,t_2]\backslash E_{x,h,S_1}$. Thus
\bean \int_x^{t_2} s(t)dt &\le& S_1\mu([x,t_2]\cap E_{x,h,S_1})-S_1 \mu([x,t_2]\backslash E_{x,h,S_1}) \\
  &=& -S_1(t_2-x) + 2S_1\mu([x,t_2]\cap E_{x,h,S_1}) \\
  &=& S_1(x-t_2+2\mu([x,t_2]\cap E_{x,h,S_1}))
\eean
In view of (\ref{eq-t2}) and $t_2>T-\ve$ (with $\ve>0$ arbitrary), we have 
$x-t_2<x-T+\ve\le 2M'/S_2-h+\ve$, thus we continue the preceding inequality as
\bean \cdots & < & S_1\left(\frac{2M'}{S_2}-h+\ve + 2\mu([x,t_2]\cap E_{x,h,S_1})\right)
\eean
By our assumption that (i) is false, we have $\mu([x,t_2]\cap E_{x,h,S_1})\le\mu(E_{x,h,S_1})<e$. 
Thus choosing $\ve$ such that $0<\ve<2(e-\mu([x,t_2]\cap E_{x,h,S_1}))$, we have
\bean \cdots &<& S_1\left(\frac{2M'}{S_2}-h+2e\right). \eean
Combining this with (\ref{eq-h}), we finally obtain $\int_x^{t_2} s(t)dt<-2M'$, which contradicts
the assumption (\ref{eq-s2}). Thus there is a point $t_1\in[x,t_2]$ such that $s(t_1)\ge S_1$. In
view of $s(t_1)\ge S_1>S_2\ge s(t_2)$, we have  $t_1\ne t_2$, thus $t_1<t_2$.
\qed

\noindent{\it Proof of Proposition \ref{prop-E}.}  Assuming that $S=\limsup |s(x)|>0$, choose
$S_1, S_2$ such that $0<S_2<S_1<S$. Then $e:=\log\frac{S_1+M}{S_2+M}>0$. Let $h$ satisfy
(\ref{eq-h}). Assume that there is an $x\ge 0$ such that $\mu(E_{x,h,S_1})<e$.  Then Lemma
\ref{lem-E1} implies the existence of $t_1,t_2$ such that $x\le t_1<t_2\le x+h$ and
$s(t_1)\ge S_1, \ s(t_2)\le S_2$. But then Corollary \ref{coro-E2} gives
$\mu([t_1,t_2]\cap s^{-1}([S_2,S_1])\ge \log\frac{S_1+M}{S_2+M}$. Since 
$[t_1,t_2]\cap s^{-1}([S_2,S_1])\subset E_{x,h,S_1}$, we have
$\mu(E_{x,h,S_1})\ge\log\frac{S_1+M}{S_2+M}=e$, which is a contradiction.  
\qed

\brem The author did not succeed in making full sense of the proof in \cite{karamata} corresponding
to that of Corollary \ref{coro-E2}. It seems that there is a logical mistake in the reasoning,
which is why we resorted to the above more topological approach.
\erem

\appendix
\section{The Prime Number Theorem}\label{sec-pnt}
\bprop \label{prop-psi} Defining $\psi(x):=\sum_{n\le x}\Lambda(n)$, we have $\psi(x)\sim x$.
\eprop

\prf Since $\Lambda(x)\ge 0$, we have that $\psi$ is non-negative and non-decreasing. Furthermore,
\be \sum_{n\le x}\psi\left(\frac{x}{n}\right)=\sum_{n\le x}\sum_{m\le x/n}\Lambda(m)
  =\sum_{r\le x}\sum_{s|r}\Lambda(s)=\sum_{r\le x}\log  r =x\log x-x+O(\log x) \label{eq-ch}\ee
by Proposition \ref{prop-Lambda1}(i) and (\ref{s2}). Now Theorem \ref{theor-LI}
implies $\psi(x)=x+o(x)$, or $\psi(x)\sim x$. 
\qed

Note that we still used only (i) of Proposition \ref{prop-Lambda1}, but now we will need (iii):

\btheor Let $\pi(x)$ be the number of primes $\le x$ and $p_n$ the $n$-th prime. Then
\bean \pi(x)&\sim&\frac{x}{\log x}, \\
 p_n &\sim &n\log n. \eean
\etheor

\prf Using Proposition \ref{prop-Lambda1}(iii), we compute
\[  \psi(x)=\sum_{n\le x}\Lambda(n)=\sum_{p^k\le x} \log p
   =\sum_{p\le x} \log p \left\lfloor\frac{\log x}{\log p}\right\rfloor \le \pi(x)\log x.  \]
If $1<y<x$ then
\[ \pi(x)-\pi(y)=\sum_{y<p\le x}1\le\sum_{y<p\le x}\frac{\log p}{\log y}\le \frac{\psi(x)}{\log y}. \]
Thus $\pi(x)\le y+\psi(x)/\log y$. Taking $y=x/\log^2x$ this gives 
\[\label{ineq} \frac{\psi(x)}{x}\le\frac{\pi(x)\log x}{x}
     \le\frac{\psi(x)}{x}\,\frac{\log x}{\log(x/\log^2x)}+\frac{1}{\log x},\]
thus $\psi(x)\sim\pi(x)\log x$. Together with Proposition \ref{prop-psi}, this gives $\pi(x)\sim x/\log x$.

Taking logarithms of $\pi(x)\sim x/\log x$, we have $\log\pi(x)\sim\log x-\log\log x\sim\log x$
and thus $\pi(x)\log\pi(x)\sim x$. Taking $x=p_n$ and using $\pi(p_n)=n$ gives $n\log n\sim p_n$. 
\qed

\brem  Karamata's proof of the Landau-Ingham theorem obviously is modeled on Selberg's original
elementary proof \cite{selberg} of the prime number theorem. However, Selberg worked with $f=\psi$
from the beginning. Most later proofs follow Selberg's approach, but there are some that work with
$M$ instead of $\psi$. Cf.\ the papers \cite{postn,kalecki} and the textbooks \cite{GL, ellison}. 
As mentioned in the introduction, the result for $M$ also follows easily from Theorem \ref{theor-LI}:
\erem

\bprop Defining $M(x)=\sum_{n\le x}\mu(n)$, we have $M(x)=o(x)$.
\eprop

\prf We define $f(x)=M(x)+\lfloor x\rfloor$, which is non-negative and non-decreasing. Now
\bean F(x) &=& \sum_{n\le x} M\left(\frac{x}{n}\right) +\left\lfloor\frac{x}{n}\right\rfloor
   =\sum_{n\le x}\sum_{m\le x/n}(\mu(m)+1) \\
   &=& \sum_{m\le x}(\mu(m)+1)\left\lfloor\frac{x}{m}\right\rfloor
   =1+\sum_{m\le x}\left\lfloor\frac{x}{m}\right\rfloor,
\eean
where the last identity is just the first in (\ref{eq-mu}). The remaining sum is known from
Dirichlet's divisor problem and can be computed in elementary fashion, 
\be \sum_{m\le x}\left\lfloor\frac{x}{m}\right\rfloor= x\log x+(2\gamma-1)x+O(\sqrt{x}), \label{eq-dir}\ee
cf.\ e.g.\ \cite{TMF}. Thus $F(x)=x\log x+(2\gamma-1)x+O(\sqrt{x})$, and Theorem \ref{theor-LI}
implies $f(x)=x+o(x)$, thus $M(x)=o(x)$.  
\qed

\brem Note that we had to define $f(x)=M(x)+\lfloor x\rfloor$ and use (\ref{eq-dir}) since
$f(x)=M(x)+x$ is non-negative, but not non-decreasing. One can generalize Theorem \ref{theor-LI}
somewhat so that it applies to functions like $f(x)=M(x)+x$ weakly violating monotonicity. 
%(As to Theorem \ref{theor-2}, this was noted in Remark \ref{rem-th2}.3.) 
But the additional effort would exceed that for the easy proof of (\ref{eq-dir}). 
\erem


\begin{thebibliography}{99}
\bibitem{balog} A. Balog: An elementary Tauberian theorem and the prime number theorem. Acta
  Math. Acad. Sci. Hung. {\bf 37}, 285--299 (1981).
\bibitem{BGT} N. H. Bingham, C. M. Goldie, J. L. Teugels: {\it Regular variation}. Cambridge
  University Press, 1987.
\bibitem{diamond} H. G. Diamond: Elementary methods in the study of the distribution of prime
  numbers. Bull. Amer. Math. Soc.  {\bf 7}, 553--589 (1982).
\bibitem{ellison} W. \& F. Ellison: {\it Prime numbers}. Wiley, Hermann, 1985. (Original French
  version: W. Ellison, M. Mend\`es France: {\it Les nombres premiers}. Hermann, 1975.)
\bibitem{erdos} P. Erd\"os: On a new method in elementary number theory which leads to an elementary proof of the
  prime-number theorem. P.N.A.S. {\bf 35}, 374--384 (1949).
\bibitem{erdos2} P. Erd\"os: On a Tauberian theorem connected with the new proof of the prime number
  theorem. J. Indian Math. Soc. (N.S.) {\bf 13}, 131--144 (1949). 
\bibitem{EI} P. Erd\"os, A. E. Ingham: Arithmetical Tauberian theorems. Acta Arith. {\bf 9}, 341--356 (1964).
\bibitem{EK} P. Erd\"os, J. Karamata: Sur la majorabilit\'e $C$ des suites des nombres
  r\'eels. Publ. Inst. Math. Acad. Serbe {\bf 10}, 37--52 (1956).
\bibitem{GL} A. O. Gelfond, Yu. V. Linnik: {\it Elementary methods in the analytic theory of
  numbers}. Pergamon Press, 1966.
\bibitem{gordon} B. Gordon: On a Tauberian theorem of Landau. Proc. Amer. Math. Soc.  {\bf 9},
  693--696 (1958). 
\bibitem{hewitt} E. Hewitt: Integration by parts for Stieltjes integrals. Amer. Math. Monthly 
  {\bf 67}, 419--423 (1960).
\bibitem{HS} E. Hewitt, K. Stromberg: {\it Real and abstract analysis}. Springer, 1965.
\bibitem{HT} A. Hildebrand, G. Tenenbaum: On some Tauberian theorems related to the prime number
  theorem. Compos. Math. {\bf 90}, 315--349 (1994).
\bibitem{ingham} A. E. Ingham: Some Tauberian theorems connected with the prime number
  theorem. 
J. London Math. Soc. {\bf 20},  171--180 (1945).
\bibitem{kalecki} M. Kalecki: A simple elementary proof of $M(x)=\sum_{n\le x}\mu(n)=o(x)$. Acta
  Arithm. {\bf 13}, 1--7 (1967).
\bibitem{karamata} J. Karamata: Sur les inversions asymptotiques de certains produits de
  convolution. Bull. Acad. Serbe Sci. (N.S.) Cl. Sci. Math.-Nat. Sci. Math. {\bf 3}, 11--32 (1957).
\bibitem{karamata2} J. Karamata: Sur les proc\'ed\'es de sommation intervenant dans la th\'eorie des
  nombres. pp. 12--31. {\it Colloque sur la th\'eorie des suites tenu \`a Bruxelles du 18 au 20 d\'ecembre
  1957}. Gauthier-Villars, 1958.
\bibitem{korevaar} J. Korevaar: {\it Tauberian theory. A century of developments}. Springer, 2004.
\bibitem{koukou} D. Koukoulopoulos: Pretentious multiplicative functions and the prime number
  theorem for arithmetic progressions. Compos. Math. {\bf 149}, 1129--1149 (2013).
\bibitem{landau} E. Landau: {\it Handbuch der Lehre von der Verteilung der Primzahlen}. Teubner,
  1909. 
\bibitem{levinson} N. Levinson: A motivated account of an elementary proof of the prime-number
  theorem. Amer. Math. Monthly {\bf 76}, 225--245 (1969).
\bibitem{nev} V. Nevanlinna: {\it \"Uber die elementaren Beweise der Primzahls\"atze und deren
    \"aquivalente Fassungen}. Univ. of Helsinki Thesis, 1964.
\bibitem{nik} A. Nikoli\'c: The story of majorizability as Karamata's condition of convergence for
  Abel summable series. Hist. Math. {\bf 36}, 405--419 (2009).
\bibitem{pollack} P. Pollack: {\it Not always buried deep. A Second Course in Elementary Number
    Theory}. Amer. Math. Soc., 2009.
\bibitem{postn} A. G. Postnikov, N. P. Romanov: A simplification of Selberg's elementary proof of the
  asymptotic law of distribution of primes. (In Russian.) Uspek. Mat. Nauk. 10:4 (66) 75--87 (1955).
\bibitem{schwarz} W. Schwarz: {\it Einf\"uhrung in Methoden und Ergebnisse der
    Primzahltheorie}. Bibliographisches Institut, 1969.
\bibitem{selberg} A. Selberg: An elementary proof of the prime-number theorem. Ann. Math. {\bf 50}, 305--313 (1949).
\bibitem{TI} T. Tatuzawa, K. Iseki: On Selberg's elementary proof of the prime number
  theorem. Proc. Japan Acad. {\bf 27}, 340--342 (1951).
\bibitem{TMF}  G. Tenenbaum, M. Mend\`es France: {\it The prime numbers and their distribution}. Amer. Math. Soc.,
  2000. (Translation of {\it Les nombres premiers}, Presses Universitaires de France, 1997.)
\bibitem{zagier} D. Zagier: Newman's short proof of the prime number theorem. 
Amer. Math. Monthly {\bf 104}, 705--708 (1997). 
\end{thebibliography}
\end{document}